\documentclass[reqno,twoside]{amsart}
\usepackage{amsmath}
\usepackage{amsfonts}
\usepackage{amssymb}
\usepackage{enumerate}

\usepackage{a4wide,amsmath,amssymb,latexsym,amsthm}

\usepackage{mathrsfs,mathtools,epic,bm}
\usepackage{hyperref}
\hypersetup{colorlinks=true,linkcolor=blue,filecolor=mangeta,urlcolor= cyan}

 \usepackage[pagewise]{lineno}

\newtheorem{theorem}{Theorem}[section]
\newtheorem{proposition}[theorem]{Proposition}
\newtheorem{lemma}[theorem]{Lemma}
\newtheorem{corollary}[theorem]{Corollary}
\theoremstyle{definition}
\newtheorem{definition}[theorem]{Definition}

\newtheorem{remark}[theorem]{Remark}

\numberwithin{equation}{section}

\def \dis {\displaystyle}

\def \R {\mathbb{R}}

\def \O {\mathcal{O}}
\def \V {\mathbb{ V}}

\def \U {\mathcal{U}}

\def \dx{\mathrm{d}x}
\def \dt{\mathrm{d}t}
\def \dy {\mathrm{d}y}
\def \dq {\mathrm{d}x \, \mathrm{d}t}

\def \ro {\rho}

\def \O {\mathcal{O}}

\def \hvarphi \widehat{\varphi}

\def\RR{{\mathbb{R}}}

\def\Om{\Omega}

\keywords{Fractional Laplacian, bilinear system, optimal control, optimality system, maximum principle, first and second order optimality conditions}


\begin{document}	
\title{Bilinear optimal control for a fractional diffusive equation}

\author{Gis\`{e}le Mophou}
\address{Gis\`ele Mophou, Laboratoire L.A.M.I.A., D\'epartement de Math\'ematiques et Informatique, Universit\'{e} des Antilles, Campus Fouillole, 97159 Pointe-\`a-Pitre,(FWI), Guadeloupe, Laboratoire  MAINEGE, Universit\'e Ouaga 3S, 06 BP 10347 Ouagadougou 06, Burkina Faso.}
 \email{gisele.mophou@univ-antilles.fr}

 	\author{Cyrille Kenne}
 \address{Cyrille Kenne,  Laboratoire L.A.M.I.A., D\'epartement de Math\'ematiques et Informatique, Universit\'{e} des Antilles, Campus Fouillole, 97159 Pointe-\`a-Pitre,(FWI), Guadeloupe, Laboratoire L3MA, DSI et IUT, Universit\'e des Antilles, Schoelcher, Martinique.}
 \email{kenne853@gmail.com}

\author{Mahamadi Warma}
\address{Mahamadi Warma, Department of Mathematical Sciences and the Center for Mathematics and Artificial Intelligence (CMAI), George Mason University,  Fairfax, VA 22030, USA.}
 \email{mwarma@gmu.edu}

\thanks{The third  author is partially supported by  US Army Research Office (ARO) under Award NO: W911NF-20-1-0115}

\begin{abstract}
	We consider a bilinear optimal control for an evolution equation involving the fractional Laplace operator of order $0<s<1$. We first give some existence and uniqueness results for the considered evolution equation. Next, we establish some weak maximum principle results allowing us to obtain more regularity of our state equation.  Then, we consider an optimal control problem which consists to bring the state of the system at final time to a desired state. We show that this optimal control problem has a solution and we derive the first and second order optimality conditions. 
\end{abstract}
	\maketitle

	\section{Introduction}
	Let $\Om\subset\R^N$ ($N\ge 1$) be a bounded domain with  boundary $\partial \Om$ and $\omega\subset\Omega$ a nonempty open set. 	Given  $T>0$, $ \alpha>0$ and  $\rho^d\in L^\infty(\Omega),$ we  are interested  to the optimal control problem: Find
	\begin{equation}\label{pbopt}
	\inf_{v\in \U} J(v):=\frac{1}{2}\|\rho(\cdot,T)-\rho^d\|^2_{L^2(\Omega)}+ \frac{\alpha}{2}\|v\|^2_{L^2(\omega\times(0,T))},
	\end{equation}
subject to the constraints that $\rho$ solves the
space fractional diffusion equation
	\begin{equation}\label{modelc}
	\left\{
	\begin{array}{rllll}
	\dis  \rho_t+(-\Delta)^s\rho &=& v\rho\chi_{\omega } \qquad &\hbox{in}& Q:= \Omega\times (0,T),\\
	\rho&=&0  &\hbox{in}& \Sigma:=(\R^N\setminus \Omega)\times (0,T) ,\\
	\rho(\cdot,0)&=& \rho^0 &\hbox{in}& \Omega,
	\end{array}
	\right.
	\end{equation}
 and the set of admissible controls is given by	
\begin{equation}\label{defuad}
	\U:=\left\{v\in L^\infty(\omega\times(0,T)): m\leq v\leq M\right\},
	\end{equation}
	where $ m,M\in \R,\, M>m$.
In \eqref{modelc}, $(-\Delta)^s$ denotes the fractional Laplace operator of order $0<s<1$, $\rho^0\in L^\infty(\Omega)$ and  $\chi_\omega$ is the characteristic function of  $\omega$.\par	
Bilinear systems are used to describe many processes in biology, ecology and engineering. These systems which are nonlinear due to the product between the input and the state variable are gained in interest to many researchers. We refer to Bruni et \textit{al.} \cite{bruni1974} for instance. Optimal control of such systems has been widely  investigated.  In the case of control depending only on time,  Bradley et \textit{al.} \cite{lenhart1999} proved the existence and uniqueness of a bilinear optimal control. Actually, the control which acts as a multiplier of a velocity term is a positive uniformly bounded function of time. The uniqueness of the optimal control were achieved for a time $T$ sufficiently small. Addou et \textit{al.} \cite{addou2002} studied a bilinear optimal control of a system governed by  a fourth-order parabolic operator. The bilinearity appeared in the form of the scalar product of the vector of controls and the gradient of the state. The authors proved under suitable hypotheses the existence of an optimal control that they characterized with an optimality system. Then, assuming that the initial state is small enough, they obtained the uniqueness of the optimal control. In \cite{ito2007}, an optimal bilinear control of an abstract Schr\"odinger equation was considered. The existence of an optimal control depending only on time is proved and the first order optimality system is derived. The paper \cite{zerrik2014} considered a regional quadratic control problem for distributed bilinear systems. They proved that an optimal control exists and  gave an optimality system that characterizes the control. Considering a control depending on time and space,   \cite{lenhart1993} studied an optimal control of linear heat equation with convective boundary condition in which  the heat transfer is took as the control. They proved the existence and uniqueness of the optimal control  and the solution of the optimality system. The results were achieved by means of compactness  and  maximum principle results.  In \cite{zerrik2020} they investigated a constrained regional control problem of a bilinear plate equation. They proved using some compactness results the existence of an optimal control that they characterized with an optimality system. The cases of time or space dependent control were also discussed.
Recently,  \cite{aronna2021} studied an optimal control problem subject to the Fokker-Planck equation. They proved the existence of optimal controls and derived the first and second order optimality conditions.
We refer to \cite{aronna2018, lenhart1999, lingfei2021, sidi2013, ouzara2015, shen2013, zerrik2019, zerrik2016, ztot2011} for more details on bilinear optimal control problems involving PDEs of integer order,  and to \cite{bahaa2016, bahaa2017, bahaa2018, mophou2011, mophou2011a} and the references therein for additional information on optimal control problems involving fractional-order operators.\par
Actually, since bilinear controlled PDEs are nonlinear, the most challenging issue in controlling such models is not only to find appropriate compactness results to obtain the existence of an optimal control, but also to derive necessary and sufficient optimality conditions and  the uniqueness of optimal solutions. This can be achieved by improving the regularity of solutions to the optimality systems.\par
To be the best of our knowledge this is the first work on control problems associated to fractional bilinear PDEs.  We have obtained the following specific results:
	\begin{itemize}
\item Corollary \ref{corexistence} and Theorem \ref{theominmax2} show existence, uniqueness and regularity of solutions of  \eqref{modelc}. The regularity  is obtained by proving some results of maximum principle that are interesting in their own independently of the application given here.
\item Theorem \ref{existcontrol1} gives the existence of solutions to the control problem \eqref{modelc}-\eqref{pbopt}.
\item The first order necessary optimality conditions are given in Theorem \ref{theoSO}.
\item The second order necessary and sufficient conditions are contained in Theorems \ref{theo-514} and  \ref{Quadratic growth}.
\end{itemize}

Let us observe that control problems of elliptic or parabolic PDEs associated with the fractional Laplace operator are always delicate due to the lack of known enough regularity results for solutions of the associated state and adjoint state equations. That are important to analyze the control problem.  To be more specific even on smooth domains,  the eigenfunctions of the operator $(-\Delta)_D^s$  defined in \eqref{FDL} below (the realization of the fractional Laplace  operator with the zero Dirichlet exterior condition)  only enjoy the following known results.  They belong to $W^{2s,2}(\Omega)$ if $0<s\le 1/2$ and they belong to $W^{\frac 12+s-\varepsilon,2}(\Omega),$ for every $\varepsilon>0$,  if $1/2\le s<1$.  We refer to \cite{Grub,G-JFA} for more information.  For the associated parabolic problem we refer to \cite{BWZ3-1} and the references therein. 
For the local case $s=1$, it is well-known that eigenfunctions are smooths.
In the present paper we have proved some space-time regularity results for the state  equation \eqref{modelc} and the associated adjoint equation. Also a maximum principle has been obtained. These results are crucial to analyze the existence of optimal solutions and to characterize the optimality conditions.

The paper is organized as follows. In Section \ref{prelim}, we give some definitions and properties of the fractional Laplacian and some known results. In Section \ref{existe}, we first prove the existence, regularity and uniqueness  of weak solutions  to the bilinear fractional diffusive system. Then, under some assumptions on the data, we establish appropriate maximum principle results. We prove in Section \ref{control} that there exists at least one optimal control solution of \eqref{modelc}-\eqref{pbopt}. In Section \ref{sec-5} we derive the first and second order optimality conditions and systems.  
\section{Preliminaries} \label{prelim}	
For the sake of completeness,  we give some well-known results that are used throughout the paper.
	We start by introducing the fractional Laplace operator. Given $0<s<1$,  we let
	\begin{align*}
	\mathcal L_s^{1}(\RR^N):= \left\{w:\RR^N\to\RR\;\hbox{ measurable and}\; \int_{\RR^N}\frac{|w(x)|}{(1+|x|)^{N+2s}}\;\dx<\infty\right\}.
	\end{align*}
	For $w\in \mathcal L_s^{1}(\RR^N)$ and $\varepsilon>0$, we set
	\begin{align*}
	(-\Delta)_\varepsilon^s w(x):= C_{N,s}\int_{\{y\in\RR^N:\;|x-y|>\varepsilon\}}\frac{w(x)- w(y)}{|x-y|^{N+2s}}\;\dy,\;\;x\in\RR^N,
	\end{align*}
	where $C_{N,s}$ is a normalization constant given by
	$\dis C_{N,s}:= \frac{s2^{2s}\Gamma\left(\frac{2s+N}{2}\right)}{\pi^{\frac{N}{2}}\Gamma(1-s)}$.
	The fractional Laplacian $(-\Delta)^s$ is defined by the following singular integral:
	\begin{align}\label{fl_def}
	(-\Delta)^sw(x):= C_{N,s}\,\hbox{P.V.}\int_{\RR^N}\frac{w(x)-w(y)}{|x-y|^{N+2s}}\;\dy =
	\lim_{\varepsilon\downarrow 0}(-\Delta)_\varepsilon^s w(x),\;\;x\in\RR^N,
	\end{align}
	provided that the limit exists for a.e. $x\in\RR^N$.
We refer to \cite{nezza2012,kwasnicki2017} for  equivalent definitions of $(-\Delta)^s$.\par
Next, we introduce the function spaces needed to investigate our problem.
Let $\Omega\subset\RR^N$ ($N\ge 1$) be an arbitrary open set and $0<s<1$. We define the fractional order Sobolev space
\begin{align*}
H^s(\Omega):=\left\{u\in L^2(\Omega):\; \int_{\Omega}\int_{\Omega}\frac{|u(x)-u(y)|^2}{|x-y|^{N+2s}}\;dxdy<\infty\right\}
\end{align*}
and we endow it with the norm given by
\begin{align*}
\|u\|_{H^s(\Omega)}=\left(\int_{\Omega}|u|^2\;dx+\int_{\Omega}\int_{\Omega}\frac{|u(x)-u(y)|^2}{|x-y|^{N+2s}}\;dxdy\right)^{1/2}.
\end{align*}
We set
	$$
	H_0^{s}(\Omega):=\Big\{w\in H^s(\R^N):\;w=0\;\hbox{ in }\;\RR^N\setminus\Omega\Big\}.
	$$
	Then,  $H_0^{s}(\Omega)$ endowed with the norm
	\begin{equation}\label{norm0bisinter}
	\|w\|_{H_0^{s}(\Omega)}:=\left(\frac{C_{N,s}}{2}\int_{\R^N}\int_{\R^N}
	\frac{(w(x)-w(y))^2}{|x-y|^{N+2s}}\;\dx\, \dy\right)^{1/2},
	\end{equation}
	is a Hilbert space (see e.g. \cite[Lemma 7]{servadei}). We let $H^{-s}(\Omega):=(H_0^s(\Omega))^\star$ be the dual space of $H_0^s(\Omega)$ with respect to the pivot space $L^2(\Omega)$, so that we have the following continuous embeddings (see e.g. \cite{ATW}):
	\begin{equation}\label{injection1}
	H_0^{s}(\Omega)\hookrightarrow L^2(\Omega)\hookrightarrow H^{-s}(\Omega).
	\end{equation}	
	From now on,  for any $\rho,\psi\in H_0^{s}(\Omega)$, we set
\begin{equation}\label{defF}
	\mathcal{F}(\rho,\psi):= \frac{C_{N,s}}{2} \int_{\R^N}\int_{\R^N}
	\frac{(\rho(x)-\rho(y))(\psi(x)-\psi(y))}{|x-y|^{N+2s}}\;\dx\, \dy.
	\end{equation}
	Hence, the norm on $H_0^{s}(\Omega)$ given by \eqref{norm0bisinter} becomes
	$\|w\|_{H_0^{s}(\Omega)}=\left(\mathcal{F}(w,w)\right)^{1/2}$.\par
We let the operator $(-\Delta)_D^s$ on $L^2(\Omega)$ be given by
\begin{equation}\label{FDL}
D((-\Delta)_D^s):=\{u\in H_0^s(\Omega):\;(-\Delta)^su\in L^2(\Omega)\},\; (-\Delta)_D^su:=(-\Delta)^su\;\hbox{ in }\Omega.
\end{equation}
Then,  $(-\Delta)_D^s$ is the realization in $L^2(\Omega)$ of $(-\Delta)^s$ with the zero Dirichlet exterior condition.\par
The following result is well-known (see e.g.  \cite{Cl-Wa,GW-CPDE}).
\begin{proposition}\label{Prop-22}
Let $(-\Delta)_D^s$ be the operator defined in \eqref{FDL}.
Then, $(-\Delta)_D^s$ can be also viewed as a bounded operator from $H_0^s(\Omega)$ into $H^{-s}(\Omega)$ given by
\begin{align}\label{FDLB}
\langle (-\Delta)_D^su,v\rangle_{H^{-s}(\Omega),H_0^s(\Omega)}=\mathcal F(u,v),\;\;\; u,v\in H_0^s(\Omega).
\end{align}
\end{proposition}	
We also need the following compactness result.	
	\begin{theorem}\label{thmcompact}\cite[Theorem 5.1, Page 58]{lions1969}
		Let $B_0,B,B_1$ be three Banach spaces such that we have the continuous embeddings $B_0\hookrightarrow B\hookrightarrow B_1$, with $B_i$ being reflexive, $i=0,1$. Assume that the embedding $B_0\hookrightarrow B$ is also compact and set
		$$\mathcal{W}:=\dis \left\{\rho \in L^2((0,T);B_0):\;\dis  \rho_t\in L^2\left((0,T);B_1\right)\right\},$$
		with $T<\infty$. Then, $\mathcal{W}$ is compactly embedded in $L^2((0,T);B)$.
	\end{theorem}	
	Let us recall the following result given in \cite[Page 37]{lions2013}.	
	\begin{theorem}\label{Theolions61} 
		Let $\left(F, \|\cdot\|_F\right)$ be a Hilbert space. Let $\Phi$ be a subspace of $F$ endowed with a pre-Hilbert scalar product $(((\cdot,\cdot)))$, with associated norm $|||\cdot|||$ .  Moreover, let $\mathcal E:F\times \Phi\to \mathbb{C}$ be a sesquilinear form.  Assume that the following hypotheses hold:
		\begin{enumerate}
			\item  The embedding $\Phi \hookrightarrow F$ is continuous, i.e., there is a constant $C_1>0$ such that
			\begin{equation*}
			\|\varphi\|_{F}\leq C_1|||\varphi|||~~~\hbox{for all}~ \varphi~ \in~ \Phi.
			\end{equation*}		
			\item For all $\varphi\in \Phi$, the mapping $u\mapsto \mathcal E(u,\varphi)$ is continuous on $F$.		
			\item There is a constant $C_2>0$  such that
			\begin{equation*}
			\mathcal{E(\varphi,\varphi)}\geq C_2 |||\varphi|||^2~~~\hbox{for all}~~\varphi\in \Phi.
			\end{equation*}
		\end{enumerate}
		If $\varphi\mapsto L(\varphi)$ is a  continuous linear functional  on $\Phi$,  then there exists $u\in F$  verifying
		$\mathcal	E(u,\varphi)=L(\varphi)~~ \text{for all} ~~\varphi\in \Phi.$ In addition, if $(\Phi,\||\cdot\||)$ is a Hilbert space,  then the solution is unique.
	\end{theorem}	
	Next, let $\mathbb{X}$ be a Banach space with dual $\mathbb{X}^\star$. We set
	\begin{equation}\label{defW0T}
	W(0,T;\mathbb{X}):= \left\{\psi \in L^2((0,T);\mathbb{X}): \psi_{t} \in L^2\left((0,T);\mathbb{X}^\star\right)\right\}.
	\end{equation}
	Then $W(0,T;\mathbb{X})$ endowed with the  norm given by
	\begin{equation}\label{normW0T}
	\|\psi\|^2_{W(0,T;\mathbb{X})}=\|\psi\|^2_{L^2(0,T;\mathbb{X})}+\|\psi_t\|^2_{
		L^2\left(0,T;\mathbb{X}^\star\right)},
	\end{equation}
	is a Hilbert space.  Moreover,  if $\mathbb{Y}$ is a Hilbert space that can be identified with its dual $\mathbb{Y}^\star$  and we have the continuous embeddings $\mathbb{X}\hookrightarrow \mathbb{Y}=\mathbb Y^\star\hookrightarrow \mathbb{X}^\star,$
	then using \cite[Theorem 1.1, page 102]{lions1971},  we have the continuous embedding
	\begin{equation}\label{contWTA}
	W(0,T;\mathbb{X})\hookrightarrow C([0,T];\mathbb{Y}).
	\end{equation}
	\section{Existence results and maximum principle}\label{existe}	
Throughout the rest of the paper $\Omega\subset\R^N$ is an arbitrary bounded domain.  For a nonempy open set $\omega\subset\Omega$ and $T>0$, we shall let $\omega_T:=\omega\times (0,T)$ and we shall simply denote by $\|\cdot\|_{\infty}$ the $L^\infty$-norm in $L^\infty(\omega_T)$.  In addition, for $u,w\in H_0^s(\Omega)$, we let $\mathcal F(u,w)$ denote the bilinear form given in \eqref{defF}.  To symplify the notations, we set \begin{equation}
	\mathbb{V}:=  H^s_0(\Omega) \hbox{ and }\mathbb V^\star:=H^{-s}(\Omega),
	\end{equation}
and we let $ \left\langle\cdot,\cdot\right\rangle_{\V^\star,\V}$ denote the duality mapping between $\V^\star$ and $\V$.
	\subsection{Existence results}\label{existence}
	For $r>0$ a real number, we consider the system
	\begin{equation}\label{modelint}
	\left\{
	\begin{array}{rllll}
	\dis  z_t+(-\Delta)^sz +rz&=& vz\chi_{\omega }+e^{-rt}f \qquad &\hbox{in}& Q,\\
	z&=&0  &\hbox{in}& \Sigma ,\\
	z(\cdot,0)&=& \rho^0 &\hbox{in}& \Omega.
	\end{array}
	\right.
	\end{equation}
	\begin{definition}\label{weaksolutionint}
		Let $f\in L^2((0,T);\V^\star))$, $v\in L^\infty(\omega_T)$ and $\rho^{0}\in L^2(\Omega)$. We  say that
		$z\in L^2((0,T);\V)$ is a weak solution of \eqref{modelint}, if the  equality
		\begin{equation}\label{Eq-Def31int}
\begin{array}{ccccc}
	\dis -\int_0^T \langle \phi_t,z\rangle_{\V^\star,\V}\, dt &+ \dis\int_0^T \mathcal{F}(z,\phi)dt+r\int_{Q}z\phi\, \dq-\int_{\omega_T}vz\, \phi \;\dq\notag\\
		&=\dis\int_0^Te^{-rt}\langle f,\, \phi\rangle_{\V^\star,\V}\; dt+\int_\Omega \rho^0\,\phi(0) \dx,
\end{array}
\end{equation}
	holds, for every $\phi \in H(Q)$,  where
	\begin{equation}\label{FP}
H(Q):=\left\{\varphi\in W(0,T;\V) \hbox{ and } \varphi(\cdot,T)=0 \hbox{ a.e.  in } \Omega\right\}.
	\end{equation}
	\end{definition}
We have the following existence result.	
	\begin{theorem}\label{theoremexistenceint}
		Let $f\in L^2((0,T);\V^\star)$, $v\in L^\infty(\omega_T)$,  $r=\|v\|_{\infty}$ and $\rho^{0}\in L^2(\Omega)$.
Then, there exists a unique weak solution
$z \in W(0,T;\V)$ of  \eqref{modelint}.
		In addition, there is a constant $C=C(N,s,\Omega)>0$  such that
\begin{subequations}\label{estimationintw0T}
\begin{alignat}{11}
\sup_{\tau\in [0,T]}\|z(\cdot,\tau)\|^2_{L^2(\Omega)}\leq  \|f\|^2_{L^2((0,T);\V^\star)}+\|\rho^0\|^2_{L^2(\Omega)},
\label{estimationintw0TS1}\\
\|z\|^2_{L^2((0,T);\V)}\leq  \|f\|^2_{L^2((0,T);\V^\star)}+\|\rho^0\|^2_{L^2(\Omega)},
\label{estimationint2}\\
 \|z\|_{W(0,T;\V)}\leq (C\|v\|_{\infty}+3)\left(\|f\|_{L^2((0,T);\V^\star)}+\|\rho^0\|_{L^2(\Omega)}\right)
 \label{estimationintw0TS3}.
\end{alignat}
\end{subequations}
\end{theorem}
	\begin{proof}We proceed in four steps.\par 	
 \textbf{Step 1.}  We prove existence  by using Theorem \ref{Theolions61}.
Recall that the norm on  $L^2((0,T);\V)$ is given by
		$$\|z\|^2_{L^2((0,T);\V)}=\int_0^T \|z(\cdot,t)\|^2_{\V}\, dt.$$
We consider the norm defined on  $H(Q)$ by $\||z\||^2_{H(Q)}:=\|z\|^2_{L^2((0,T);\V)}+\|z(\cdot,0)\|^2_{L^2(\Omega)}.$
It is clear that we have the continuous embedding $H(Q)\hookrightarrow L^2((0,T);\V)$.\par 		
Now, let $\varphi \in H(Q)$ and consider the bilinear form $\mathcal{E}(\cdot,\cdot):L^2((0,T);\V)\times H(Q)\to\R$ given by
		\begin{equation}\label{defCalE}
		\mathcal{E}(z,\varphi):=  -\int_0^T \langle \phi_t,z\rangle_{\V^\star,\V}\, dt+\dis\int_0^T \mathcal{F}(z,\varphi)dt+r\int_{Q}z\varphi\, \dq-\int_{\omega_T}vz\, \varphi \;\dq.
		\end{equation}
Using Cauchy-Schwarz's inequality,  we get that
	\begin{align*}
		\dis |\mathcal{E}(z,\varphi)|
		\leq  \Big(\|\varphi_{t}\|_{L^2((0,T);\V^\star)}+\left(r+\|v\|_{\infty}\right)\|\varphi\|_{L^2(Q)}+\|\varphi\|_{L^2((0,T);\V)}\Big)\|z\|_{L^2((0,T);\V)}.
\end{align*}
		Consequently, for every fixed $\varphi\in H(Q),$
		the functional  $z\mapsto \mathcal{E}(z,\varphi)$ is continuous on $L^2((0,T);\V).$\par 			
		Next,  since $r=\|v\|_{\infty}$,  we have that for every  $\varphi\in H(Q)$,
$$\begin{array}{rlll}
		\mathcal{E}(\varphi,\varphi)&=&\dis -\int_0^T \langle \varphi_t,\varphi\rangle_{\V^\star,\V}\, dt+\int_0^T \mathcal{F}(\varphi,\varphi)\,\dt+r\int_{Q}\varphi^2 \,\dq-\int_{\omega_T}v\varphi^2 \;\dq\\
		&\geq& \frac 12  \|\varphi(\cdot,0)\|^2_{L^2(\Omega)}+\dis \int_0^T \mathcal{F}(\varphi,\varphi)dt+(r-\|v\|_{\infty})\int_{Q}\varphi^2 \,\dq\\
		&\geq & \frac{1}{2}\||\varphi\||^2_{H(Q)}.
\end{array}
$$
Hence, $\mathcal E$ is coercive on $H(Q)$.\par
		Finally, let us consider the functional $L:H(Q)\to \R$ defined by
		$$
		L(\varphi):=\dis \int_0^Te^{-rt}\langle f,\, \phi\rangle_{\V^\star,\V}\; dt+\int_\Omega \rho^0\,\varphi(x,0)\, \dx .
		$$
Using Cauchy-Schwarz's inequality,  there is a constant $C>0$ such that
$$
		|L(\varphi)|
		\leq  C\left( \|\rho^0\|_{L^2(\Omega)}+\|f\|_{L^2((0,T);\V^\star)}\right)\|\varphi\|_{H(Q)}.
	$$
		Therefore,  $L$ is continuous and linear on $H(Q)$. It follows from Theorem \ref{Theolions61} that there exists $z\in L^2((0,T);\V)$ such that
		$\mathcal{E}(z,\varphi)= L(\varphi),\quad \forall \varphi \in H(Q)$.
		We have shown that the system \eqref{modelint} has a solution $z\in L^2((0,T);\V)$ in the sense of Definition \ref{weaksolutionint}.\par	
\noindent \textbf{Step 2.} We show that $z_t\in L^2((0,T);\V^\star)$
Notice that  \eqref{modelint} can be rewritten as the abstract Cauchy problem
\begin{equation}\label{mod-2}
\left\{
\begin{array}{rllll}
	\dis  z_t+(-\Delta)_D^sz +rz&=& vz\chi_{\omega }+e^{-rt}f \qquad &\hbox{in}& Q,\\
	z(\cdot,0)&=& \rho^0 &\hbox{in}& \Omega,
	\end{array}
	\right.
	\end{equation}
where $(-\Delta)_D^s$ is defined in \eqref{FDL}.
Since $z\in L^2((0,T);\V)$,  it follows from Proposition \ref{Prop-22} that $(-\Delta)_D^sz(\cdot,t)\in \V^\star$.  Since $v\in L^\infty(\omega_T)$, we have $zv\chi_\omega\in L^2(\omega_T)\hookrightarrow  L^2((0,T);\V^\star)$ and we can deduce that
$z_t(\cdot,t)=-(-\Delta)^s_Dz(\cdot,t) -rz(\cdot,t)+ (vz)(\cdot,t)\chi_{\omega }+e^{-rt}f(\cdot,t)\in \V^\star.$\par
If we take the duality map between \eqref{mod-2} and $\phi\in L^2((0,T);\V)$,  and use Proposition \ref{Prop-22},  we obtain
	\begin{align*}
	\dis  \left\langle z_t(t),\phi(t)\right\rangle_{\V^\star,\V} + \mathcal{F}(z(t),\phi(t))+r\int_{\Omega}z(t)\, \phi(t) \,\dx
	 =&\dis\int_{\omega}v(t)z(t)\, \phi(t) \,\dx\\
&+\dis
 e^{-rt} \langle f(t),\, \phi(t)\rangle_{\V^\star,\V}.
	\end{align*}
This implies that there is a constant $C>0$ such that
\begin{equation}\label{aa}
 |\left\langle z_t(t),\phi(t)\right\rangle_{\V^\star,\V}| \leq  \dis \left[(C\|v\|_{\infty}+1)\|z(t)\|_{\V}+\|f(t)\|_{\V^\star}\right]\|\phi(t)\|_{\V}.
\end{equation}
Integrating \eqref{aa} over $(0,T)$,  we get that there are two constants $C=C(N,s,\Omega)>0$ and $C_1:= (C\|v\|_{\infty}+1)$ such that
	\begin{equation} \label{estimationint3-1}
	\dis\int_0^T |\left\langle z_t(t),\phi(t)\right\rangle_{\V^\star,\V}| dt
\leq \left[C_1\|z\|_{L^2((0,T);\V)}+\|f\|_{L^2((0,T);\V^\star)}\right]
\|\phi\|_{L^2((0,T);\V)}.
	\end{equation}
Using \eqref{estimationint2} we get from \eqref{estimationint3-1} that
\begin{equation} \label{estimationint3}
 \|z_t\|_{L^2((0,T);\V^\star)}\leq (C\|v\|_{\infty}+2)\left(\|f\|_{L^2((0,T);\V^\star)}+\|\rho^0\|_{L^2(\Omega)}\right).
\end{equation}
Thus, $z_t\in L^2((0,T);\V^\star)$ and we have shown that $z\in W(0,T;\V)$.\par
\noindent{\bf Step 3.} We show  \eqref{estimationintw0TS1},  \eqref{estimationint2} and \eqref{estimationintw0TS3}.
If we take the duality map between  \eqref{mod-2} and $z\in W(0,T;\V)$, use Proposition \ref{Prop-22} and  Young's inequality, we get
		\begin{align*}
\dis  \frac 12\frac{d}{dt}\|z(t)\|^2_{L^2(\Omega)}+\mathcal{F}(z(t),z(t))+r\|z(t)\|^2_{L^2(\Om)}
=&  e^{-rt}\langle f(t),z(t) \rangle_{\V^\star,\V}+\int_{\omega} v(t)z^2(t)\, dx\\
		\leq&\frac{1}{2}\|f(t)\|^2_{\V^\star}+\frac{1}{2}\|z(t)\|^2_{\V}+\dis \|v\|_{\infty}\|z(t)\|^2_{L^2(\Omega)}.
		\end{align*}
		Hence,
		\begin{equation}\label{inter2}
		 \frac{d}{dt}\|z(t)\|^2_{L^2(\Omega)}+	\|z(t)\|^2_{\V}\leq \dis \|f(t)\|^2_{\V^\star},
		\end{equation}
because  $r=\|v\|_{\infty}$.
Integrating \eqref{inter2} over $(0,\tau)$, with $\tau\in [0,T]$, we get that
$$
\|z(\tau)\|^2_{L^2(\Omega)}+\int_{0}^{\tau}\|z(t)\|^2_{\V}\,dt\leq \|f\|^2_{L^2((0,T);\V^\star)}+\|\rho^0\|^2_{L^2(\Omega)},
$$
from which we deduce  \eqref{estimationintw0TS1} and \eqref{estimationint2}.
Combining  \eqref{estimationint2}-\eqref{estimationint3} and recalling  \eqref{normW0T},  we get \eqref{estimationintw0TS3}.\par
 \textbf{Step 4.} We prove  uniqueness.
		Assume that there exist $z_1$ and $z_2$ solutions to \eqref{modelint} with the same right hand side $f,\, v$ and initial datum $\rho^0$.  Set $\tilde{z}:=z_1-z_2$. Then,  $\tilde z$ satisfies
		\begin{equation}\label{pazero}
		\left\{
		\begin{array}{rllll}
		\dis  \tilde{z}_t+(-\Delta)^s\tilde{z}+r\tilde{z}&=&v\tilde{z}\chi_{\omega } \qquad &\hbox{in}& Q,\\
		\tilde{z}&=& 0 &\hbox{in}& \Sigma,\\
		\tilde{z}(\cdot,0)&=& 0 &\hbox{in}& \Omega.
		\end{array}
		\right.
		\end{equation}	
From Step 2, we have that $\tilde z$ has the regularity to be taken as a test function in \eqref{pazero}.  So,  if we take the duality map of \eqref{pazero} with $\tilde{z}$,  use Proposition \ref{Prop-22} and integrate over $(0,T)$ we obtain
$$
 \frac 12\|\tilde{z}(\cdot,T)\|^2_{L^2(\Omega)}+\int_0^T \mathcal{F}(\tilde{z},\tilde{z})\, \dt+r\| \tilde{z}\|^2_{L^2(Q)} = \dis \int_{\omega_T} v\tilde{z}^2 dxdt\leq\|v\|_{\infty}\|\tilde{z}\|^2_{L^2(Q)}.$$
Since  $r=\|v\|_{\infty}$,  we can deduce  that
		$ \frac 12\|\tilde{z}(\cdot,T)\|^2_{L^2(\Omega)}+\|\tilde{z}\|^2_{L^2((0,T);\V)}\leq 0.$
		Hence, $\tilde{z}=0$  in $\R^N\times [0,T]$. Thus, $z_1=z_2$ in $\R^N\times [0,T]$ and we have shown uniqueness.
	\end{proof}
	
Next, we consider the system
\begin{equation}\label{model}
	\left\{
	\begin{array}{rllll}
	\dis  \rho_t+(-\Delta)^s\rho &=& f+v\rho\chi_{\omega } \qquad &\hbox{in}& Q,\\
	\rho&=&0  &\hbox{in}& \Sigma ,\\
	\rho(\cdot,0)&=& \rho^0 &\hbox{in}& \Omega.
	\end{array}
	\right.
	\end{equation}
We have the following result that can be viewed as a corollary of Theorem \ref{theoremexistenceint}.
	\begin{corollary}\label{theoremexistence}
		Let $f\in L^2((0,T);\V^\star)$, $v\in L^\infty(\omega_T)$ and $\rho^{0}\in L^2(\Omega)$.  Then, there exists a unique weak solution $\ro \in W(0,T;\V)$ of \eqref{model}.
		In addition, there is a constant $C=C(N,s,\Omega)>0$ such that 	
\begin{subequations}\label{estimation22W0Tout}
\begin{alignat}{11}
\sup_{\tau\in [0,T]}\|\rho(\cdot,\tau)\|^2_{L^2(\Omega)}\leq  e^{2\|v\|_{\infty}T}\left[ \|f\|^2_{L^2((0,T);\V^\star)}+\|\rho^0\|^2_{L^2(\Omega)}\right],
\label{estimation1S1}\\
\|\rho\|^2_{L^2((0,T);\V)}\leq  e^{2\|v\|_{\infty}T}\left[ \|f\|^2_{L^2((0,T);\V^\star)}+\|\rho^0\|^2_{L^2(\Omega)}\right],
\label{estimation1}\\
 \|\rho\|_{W(0,T;\V)}\leq \left[\left( 1+C\|v\|_{\infty}\right)e^{\|v\|_{\infty}T}+1\right]\left[\|f\|_{L^2((0,T);\V^\star)}+
 \|\rho^0\|_{L^2(\Omega)}\right]
 \label{estimation22W0T}.
\end{alignat}
\end{subequations}
	\end{corollary}	
	\begin{proof}
Since $\rho:=e^{\|v\|_{\infty}t}z$ is a weak solution of \eqref{model} if and only if $z$ is a weak solution of \eqref{modelint}, we have from Theorem \ref{theoremexistenceint} that there exists a unique solution $\rho \in W(0,T;\V)$ of \eqref{model}.\par
Next,  letting  $z=e^{-\|v\|_{\infty}t}\rho$ in \eqref{estimationintw0TS1} and \eqref{estimationint2} we  respectively deduce that \eqref{estimation1S1} and \eqref{estimation1} hold true.\par
Finally,  we show \eqref{estimation22W0T}.  Let $\phi\in L^2((0,T);\V)$. If we take the duality map between \eqref{model} and $\phi(t)$,  and use  Proposition \ref{Prop-22},  we obtain for a.e $t\in (0,T)$,
$$
\dis \left\langle \rho_t(t),\phi(t)\right\rangle_{\V^\star,\V} +  \mathcal{F}(\rho(t),\phi(t))=\int_{\omega}v(t)\ro(t)\, \phi(t) \,\dx+
\dis \langle f(t),\, \phi(t)\rangle_{\V^\star,\V}.
$$
Using the Cauchy-Schwarz inequality and integrating over $(0,T)$ we can deduce that
\begin{align}\label{estimation22}
\|\rho_t\|_{L^2((0,T);\V^\star)} \leq & \left( 1+C\|v\|_{\infty}\right)e^{\|v\|_{\infty}T}\left[\|f\|_{L^2((0,T);\V^\star)} +\|\rho^0\|_{L^2(\Omega)}\right],
\end{align}
where we have also used \eqref{estimation1}. Adding \eqref{estimation1} to \eqref{estimation22}  we get \eqref{estimation22W0T}.
\end{proof}
\begin{corollary}\label{corexistence}
		Let  $v\in L^\infty(\omega_T)$ and $\rho^0\in L^2(\Omega)$.  Then, there exists a unique solution $\ro \in W (0,T;\V)$ of \eqref{modelc}.
In addition,  there is a constant $C=C(N,s,\Omega)>0$ such that 	
		\begin{equation}\label{estimation22coW0T}
	\sup_{\tau\in [0,T]}\|\rho(\cdot,\tau)\|_{L^2(\Omega)}\leq C \|\rho\|_{W(0,T;\V)}\leq
\left( 1+C\|v\|_{\infty}\right)e^{\|v\|_{\infty}T}
 \|\rho^0\|_{L^2(\Omega)}.
		\end{equation}		
	\end{corollary}
	\begin{proof}  It suffice to  apply Corollary \ref{theoremexistence} with $f\equiv 0$.
	\end{proof}
	\subsection{Maximum principle}	
We give some useful results of maximum principle.
	\begin{lemma}\label{positive}
		Let   $\rho^0\in L^2(\Omega)$  be  such that $\rho^0\geq 0$ a.e.  in $\Omega$  and  $v\in L^\infty(\omega_T)$. Then,  the weak solution $\rho$ of  \eqref{modelc} satisfies $\rho\ge 0$ a.e.  in $\R^N\times [0,T]$.
	\end{lemma}
	\begin{proof}
			We write $\rho=\rho^+-\rho^-$, where $\rho^+:=\max(\rho,0)$ and
$\rho^-:=\max(0,-\rho)$. It is sufficient to show that $\rho^-=0$
a.e.  in $ \R^N\times [0,T]$. Notice that
$$		\rho^-=0\;\; \hbox{in } \Sigma\;\hbox{ and }\;
		\rho^-(\cdot,0)=\max(0,-\rho^0)=0\;\; \hbox{a.e. in } \Omega.
		$$
Moreover, we have that  $\rho^-\in W(0,T;\mathbb{ V})$ (see e.g. \cite{War}).  We set
\begin{equation}\label{O}
\begin{cases}
\O^-:=\dis \left\{x\in \R^N:  \rho(x,t)\leq 0 \hbox{ for a.e. } t\in (0,T)\right\}, \\
\O^+:=\dis \left\{x\in \R^N: \rho(x,t)> 0 \hbox{ for a.e. } t\in (0,T)\right\}.
\end{cases}
\end{equation}
If we take the duality map between  \eqref{modelc}  and $\psi\in \V$,  and use  Proposition \ref{Prop-22},  we get  for a.e. $t\in [0,T]$,
\begin{equation}\label{SH}
\langle \rho_t (t),\psi\rangle_{\V^\star,\V}+\mathcal{F}(\rho(t),\psi)
		=\dis  \dis \int_\omega v(t)\rho(t)\psi dx.
\end{equation}
		Taking $\psi=\rho^-(\cdot,t)$ in \eqref{SH} and noticing that $\rho_t^-=\rho_t\chi_{\mathcal O^-}$ and $\rho^+\rho^-=0$, we get
		\begin{equation}\label{ajout1}
\int_{\O^-}\frac{d}{dt}|\rho^- (t)|^2\, \dx-\mathcal{F}(\rho(t),\rho^-(t))
		=\int_{\O^-} v(t)|\rho^-(t)|^2\, \dx.
		\end{equation}
Observe also that	
		\begin{equation}\label{equality}
\dis 	\mathcal{F}(\rho(t),\rho^-(t))= \mathcal{F}(\rho^+(t),\rho^-(t))-\mathcal{F}(\rho^-(t),\rho^-(t)).
		\end{equation}
It is well-known (see e.g. \cite{Cl-Wa,GW-CPDE,War}) that $\mathcal{F}(\rho^+(t),\rho^-(t))\le 0$.
Thus, it follows from \eqref{equality} that $\mathcal{F}(\rho(t),\rho^-(t))\le 0$.
Since $-\mathcal{F}(\rho(t),\rho^-(t))\ge 0 $, it follows from \eqref{ajout1} that
			\begin{equation}
		\dis \dis \frac{1}{2}\frac{d}{dt}\|\rho^-(t) \|^2_{L^2(\O^-)} \le \int_{\O^-} v(t)\left(\rho^-(t)\right)^2 dx\leq \dis \|v\|_{\infty}\|\rho^-(t) \|^2_{L^2(\O^-)}.
		\end{equation}
Using Gronwall's Lemma we can deduce that
		$$\|\rho^-(\cdot,t) \|^2_{L^2(\O^-)}\leq e^{2t\|v\|_{\infty}}\|\rho^-(\cdot,0)\|_{L^2(\O^-)}=0,$$
where we have also used that  $\rho^-(\cdot,0)=0$ a.e. in $\O^-$.
We have shown that
		$\rho^-=0$ a.e.  in $ \O^-\times [0,T]$. Thus, $\rho^-=0$  a.e.  in $\R^N\times [0,T]$. Consequently,  $\rho\geq 0$ a.e.  in $\R^N\times [0,T]$. The proof is finished.
	\end{proof}
We have the following maximum principle.
	\begin{theorem}\label{theominmax2}
		Let   $\rho^0\in L^\infty(\Omega)$   and $v\in L^\infty(\omega_T)$. Then,   the unique weak solution $\rho$ of \eqref{modelc} belongs to $W(0,T;\V)\cap L^\infty(\R^N\times(0,T))$ and
		\begin{equation}\label{rhoinfini2}
		\dis  \|\rho\|_{L^\infty(\R^N\times(0,T))}\leq e^{\|v\|_{\infty}T}\|\rho^0\|_{L^\infty(\Omega)}.
		\end{equation}
	\end{theorem}
\begin{proof}
We set $z:=e^{-\|v\|_{\infty}t}\rho$, where $\rho$ is the  solution of \eqref{modelc}.
  Then, $z\in W(0,T;\V)$.  We claim that
	\begin{equation}\label{ineqint}
	z\leq \|\rho^0\|_{L^\infty(\Omega)}\;\hbox{ a.e. in } \R^N\times [0,T].
	\end{equation}		
	We set $w:=\|\rho^0\|_{L^\infty(\Omega)}-z$.  Then,
	$w(x,0)=\|\rho^0\|_{L^\infty(\Omega)}-\rho^0(x)\geq 0$
	for a.e. $x\in \Omega.$ Moreover, $w$ satisfies
	\begin{equation}\label{modelwb}
	\left\{
	\begin{array}{rllll}
	\dis  w_t+(-\Delta)^sw+\|v\|_{\infty}w &=& vw\chi_{\omega }+(\|v\|_{\infty}-v \chi_{\omega })\|\rho^0\|_{L^\infty(\Omega)}  &\hbox{in}& Q,\\
	w&=&\|\rho^0\|_{L^\infty(\Omega)} &\hbox{in}& \Sigma ,\\
	w(\cdot,0)&=&\|\rho^0\|_{L^\infty(\Omega)}-\rho^0 &\hbox{in}& \Omega.
	\end{array}
	\right.
	\end{equation}
 To obtain our result, it  is sufficient to show that $w^-=0$ a.e.  in $\R^N\times [0,T]$. Since $\|\rho^0\|_{L^\infty(\Omega)}\geq 0$, we have that   $w^-= 0$ in $\Sigma$. It is also clear that $w^-\in W(0,T;\mathbb{ V})$.  Let $\O^-$ and $\O^+$ be as in \eqref{O} with $\rho$ replaced by $w$.
If we take the duality map between \eqref{modelwb} and $w^-$,  use the same argument as in the proof of Lemma \ref{positive},  we get that
\begin{align*}
-	&\frac 12 \int_{\O^-} \frac{d}{dt}|w^-(t)|^2 \,\dx+\mathcal F(w(t),w^-(t)) -\|v\|_{\infty}\int_{\mathcal O^-}(w^-(t))^2\;dx\\
	=&-\int_{\O^-} v(t)(w^-(t))^2\, \dx+ \|\rho^0\|_{L^\infty(\Omega)}\int_{\O^-} (\|v\|_{\infty}-v(t))w^-(t)\, \dx.
\end{align*}
We have shown that
\begin{align}\label{MJ}
&\frac 12 \int_{\O^-} \frac{d}{dt}|w^-(t)|^2 \,\dx-\mathcal F(w(t),w^-(t)) +\|v\|_{\infty}\int_{\mathcal O^-}(w^-(t))^2\;dx\notag\\
	=&\int_{\O^-} v(t)(w^-(t))^2\, \dx- \|\rho^0\|_{L^\infty(\Omega)}\int_{\O^-} (\|v\|_{\infty}-v(t))w^-(t)\, \dx.
\end{align}
As  in the proof of Lemma \ref{positive} we have that $\mathcal F(w(t),w^-(t))\le 0$.
Since
$$\dis \|\rho^0\|_{L^\infty(\Omega)}\int_\Omega (\|v\|_{\infty}-v(t))w^-(t)\, \dx\geq 0$$
and
\begin{align*}
-\mathcal F(w(t),w^-(x,t)) +\|v\|_{\infty}\int_{\mathcal O^-}(w^-(t))^2\;dx\ge 0,
\end{align*}
we can deduce from \eqref{MJ} that
\begin{equation}
\dis \frac{1}{2}\frac{d}{dt}\|w^-(\cdot,t) \|^2_{L^2(\O^-)}\leq \dis \|v\|_{\infty}\|w^-(\cdot,t) \|^2_{L^2(\O^-)}.
\end{equation}
 It  follows from  Gronwall's Lemma again that
$$\|w^-(\cdot,t) \|^2_{L^2(\O^-)}\leq e^{2t\|v\|_{\infty}}\|w^-(\cdot,0)\|_{L^2(\O^-)}=0,$$
where we have also used that $w^-(\cdot,0)=0$ a.e. in $\O^-$. Hence,
$w^-=0$ a.e.  in  $\O^-\times [0,T]$. We can deduce that $w\geq 0$ a.e.  in $\R^N\times [0,T]$.We have shown the claim.\par
Replacing $z$ by $e^{-\|v\|_{\infty}t}\rho$ in \eqref{ineqint}, we obtain \eqref{rhoinfini2}.
The proof is finished.
	\end{proof}
Now,  let $\rho$ be the solution of \eqref{modelc} and  consider the system
\begin{equation}\label{ad21bis}
	\left\{
	\begin{array}{lllll}
	\dis  -q_t+(-\Delta)^s q&=&vq\chi_{\omega } &\hbox{in} & Q,\\
	\dis q&=&0 &\hbox{in}&\Sigma,\\
	q(\cdot,T)&=&\ro(\cdot,T)-\ro^d &\hbox{in}& \Omega,
	\end{array}
	\right.
	\end{equation}
which can be viewed as the dual system associated with \eqref{modelc}.
\begin{corollary}\label{corominmax2}
	Let $\rho^0,\rho^d\in L^\infty(\Omega)$   and $v\in L^\infty(\omega_T)$.    Then,   \eqref{ad21bis} has a unique weak solution $q$ that belongs to $W(0,T;\V)\cap L^\infty(\R^N\times (0,T))$. In addition, there is a constant $C=C(N,s,\Omega)>0$ such that
\begin{equation}\label{qW0T}
\|q\|_{W(0,T;\V)}\leq\left( 2+C\|v\|_{\infty}\right)e^{\|v\|_{\infty}T}\left(
 \|\rho^0\|_{L^2(\Omega)}+\|\rho^d\|_{L^2(\Omega)}\right)
\end{equation}
and
\begin{equation}\label{qinfini}
		\dis  \|q\|_{L^\infty(\R^N\times(0,T))}\leq  e^{\|v\|_{\infty}T}\left(\|\rho^0\|_{L^\infty(\Omega)}+\|\ro^d\|_{L^\infty(\Omega)}\right).
	\end{equation}
\end{corollary}
\begin{proof} Making  the change of variable $t\mapsto T-t$,  we have that $\varphi(x,t):=q(x,T-t)$  satisfies
	\begin{equation}\label{ad21varphi}
	\left\{
	\begin{array}{lllll}
	\dis  \varphi_t+(-\Delta)^s \varphi&=&\tilde{v}\varphi\chi_{\omega } &\hbox{in} & Q,\\
	\dis \varphi&=&0 &\hbox{in}&\Sigma,\\
	\varphi(\cdot,0)&=&\rho(\tilde{v};\cdot,0)-\ro^d &\hbox{in}& \Omega,
	\end{array}
	\right.
	\end{equation}
where $\tilde{v}(x,t)=v(x,T-t)$.
Observing that $\rho=\rho(\tilde{v})$ is a solution of \eqref{modelc} with $\rho^0\in  L^\infty(\Omega)$ and $v=\tilde{v}$, we have from Theorem \ref{theominmax2} that $\rho\in W(0,T;\V)\cap L^\infty(\R^N\times [0,T])$.
Hence,   $\rho(\tilde{v};\cdot,0)-\ro^d=\rho^0-\ro^d\in L^\infty(\Omega)$ and  it follows from Corollary \ref{corexistence}  and Theorem \ref{theominmax2} that that there exists a unique $\varphi\in W(0,T;\V)\cap L^\infty(\R^N\times [0,T])$ solution of \eqref{ad21varphi}. Thanks to \eqref{rhoinfini2} and  \eqref{estimation22coW0T}, we have that
$$
 \|q\|_{L^\infty(\R^N\times(0,T))}=\|\varphi\|_{L^\infty(\R^N\times(0,T))} \leq e^{\|\tilde{v}\|_{\infty}T}\|\rho^0-\ro^d\|_{L^\infty(\Omega)}
$$
and there is a constant $C=C(N,s,\Omega)>0$ such that
$$
 \|q\|_{W(0,T;\V)}=\|\varphi\|_{W(0,T;\V)}\leq
\left( 2+C\|\tilde v\|_{\infty}\right)e^{\|\tilde v\|_{\infty}T}\left(
 \|\rho^0\|_{L^2(\Omega)}+\|\rho^d\|_{L^2(\Omega)}\right).
		$$
Hence, $q$ satisfies  \eqref{qW0T} and \eqref{qinfini}. The proof is finished.
\end{proof}
	\section{Existence of optimal solutions} \label{control}
We are concerned with the optimal control problem \eqref{modelc}-\eqref{pbopt}.
In view of Theorem \ref{theominmax2} and the embedding  \eqref{contWTA}  the cost function $J$ is well defined.
We define the control-to-state mapping
	\begin{equation}\label{ctso}
		G:L^\infty(\omega_T)\to W(0,T;\V)\cap L^\infty(\R^N\times (0,T)),\;\; v\mapsto G(v)=\rho
\end{equation}
		which associates to each $v\in L^\infty(\omega_T)$ the unique weak solution  $\rho$ of \eqref{modelc}. Then, the optimal control problem \eqref{modelc}-\eqref{pbopt} can be rewritten as
\begin{equation}\label{optimal}
	\inf_{v\in \U} J(v)=\frac{1}{2}\|G(v)(\cdot,T)-\rho^d\|^2_{L^2(\Omega)}+ \frac{\alpha}{2}\|v\|^2_{L^2(\omega_T)}.
	\end{equation}	
	\begin{theorem}\label{existcontrol1}
Let  $\alpha>0$, $v\in \mathcal{U}$ and $\rho^{0},\rho^d\in L^\infty(\Omega)$. Then,  there exists a solution  $u\in \mathcal{U}$ of  \eqref{optimal}, and hence, of \eqref{modelc}-\eqref{pbopt}.
	\end{theorem}
	\begin{proof}
Let $v^n\in  \mathcal{U}$ be a minimizing sequence  such that
		$$\lim_{n\to \infty}J(v^n)= \inf_{v\in\U}J(v). $$	
		Since $\rho^n:=G(v^n)$ is the state associated to the control $v^n$, there is a  constant $C>0$ independent of $n$ such that
		\begin{equation}\label{bound21}
		\|\rho^n(\cdot,T)\|_{L^2(\Omega)}\leq C\;\hbox{ and } 	\|v^n\|_{L^2(\omega_T)}\leq C.
		\end{equation}
Moreover,  $\rho^n$  satisfies
\begin{equation}\label{pn1}
		\left\{
		\begin{array}{lllll}
		\dis  (\rho^n)_t+(-\Delta)^s \rho^n&=&v^n\rho^n\chi_{\omega }&\hbox{in} & Q,\\
		\rho^n&=&0 &\hbox{in}&\Sigma,\\
		\rho^n(\cdot,0)&=&\ro^0 &\hbox{in}& \Omega.
		\end{array}
		\right.
		\end{equation}
	It follows from Corollary \ref{corexistence} that \eqref{pn1} has a unique  solution  $\rho^n\in W(0,T;\V)$ satisfying
		\begin{equation}\label{pnw1}
		\dis -\int_0^T \langle\phi_t,\rho^n \rangle_{\V^\star,\V}\, dt + \int_0^T \mathcal{F}(\rho^n,\phi)dt=
		\dis \int_{\omega_T}v^n\rho^n\, \phi \;\dq+\int_\Omega \rho^0\,\phi(0) \;dx
		\end{equation}
for every $\phi\in H(Q)$, where we recall that $H(Q)$ is given in \eqref{FP}.
Thanks to  \eqref{estimation22coW0T} and  $\|v^n\|_{\infty}\le \max\{|m|,|M|\}$,  we have that there is a constant $C>0$ independent of $n$ such that
		\begin{equation}\label{bound31}
		\|\rho^n\|_{W(0,T;\V)}\leq \dis C e^{CT} \|\rho^0\|_{L^2(\Omega)}.
		\end{equation}
We can deduce that
		\begin{equation}\label{bound41}
		\|v^n\rho^n\|_{L^2(\omega_T )}\leq C\|\rho^n\|_{L^2(Q)}
		\leq C\|\rho^n\|_{W(0,T;\V)}
		\leq C e^{CT} \|\rho^0\|_{L^2(\Omega)}.
		\end{equation}
From \eqref{bound21}-\eqref{bound41},  there exist  $\eta\in L^2(\Omega)$, $u\in L^2(\omega_T)$, $\beta\in L^2(\omega_T)$ and $\rho\in W(0,T;\V) $ such that (up to a subsequence if necessary), as $n\to \infty$, we have that
			\begin{equation}\label{c11}
		v^n\rightharpoonup u   \text{ weakly in } L^2(\omega_T),
		\end{equation}
		\begin{equation}\label{c21}
		\rho^n (\cdot,T)\rightharpoonup \eta  \text{ weakly in }   L^2(\Omega),
		\end{equation}
		\begin{equation}\label{c31}
		\rho^n\rightharpoonup \ro   \text{ weakly in }  W(0,T;\V),
		\end{equation}
		\begin{equation}\label{c41}
		v^n\rho^n\rightharpoonup \beta    \text{ weakly in }   L^2(\omega_T).
		\end{equation}
		Since $\U$ is a closed convex subset of $L^2(\omega_T)$, we have that $\U$ is weakly closed and so
	\begin{equation}\label{ajout5}
u\in \U.
	\end{equation}
It follows from Theorem \ref{thmcompact}  that the embedding $W(0,T);\V)\hookrightarrow L^2(Q)$ is compact.
Hence,  from \eqref{c31} we have that,  as $n\to\infty$,
		\begin{equation}\label{strconv1}
		\rho^n\to \ro  \text{ strongly in }  L^2(Q).
		\end{equation}
Taking \eqref{c11} and \eqref{strconv1} into account and using the weak-strong convergence, we get that,  as $n\to\infty$,
		\begin{equation}\label{c61}
		v^n\rho^n\rightharpoonup u\ro    \text{ weakly in } L^1(\omega_T).
		\end{equation}
The fact that $v^n\rho^n\rightharpoonup \beta    \text{ weakly in }   L^2(\omega_T)$ (by \eqref{c41}) together with the continuous embeddings $L^\infty(\omega_T)\hookrightarrow L^2(\omega_T)\hookrightarrow L^1(\omega_T)$ imply that  $v^n\rho^n\rightharpoonup \beta    \text{ weakly in }   L^1(\omega_T)$. Using \eqref{c61} 
		 and the uniqueness of the weak limit,  we can deduce that  $\beta=u\ro$.  We have shown that,  as $n\to\infty$,
		\begin{equation}\label{convB1}
		v^n\rho^n\rightharpoonup u\ro    \text{ weakly in }  L^2(\omega_T).
		\end{equation}
		Passing to the limit,  as $n\to \infty$,  in \eqref{pnw1}, while using \eqref{c31} and  \eqref{convB1}, we get
		\begin{equation}\label{limpnw1}
		\dis -\int_0^T\langle \phi_t,\ro\,  \dx\rangle_{\V^\star,\V}\, dt +\int_0^T \mathcal{F}(\ro,\phi)dt=
		\dis \int_{\omega_T}u\ro\, \phi \;\dq+\int_\Omega \rho^0\,\phi(0) \dx
		\end{equation}
		for every $\phi\in H(Q)$.  Thus, $\ro\in W(0,T;\V)$ is the unique solution of \eqref{modelc} with $v=u$.\par 		
		Now,  let $\phi \in W(0,T;\V)$. If we take the duality map between \eqref{pn1} and $\phi$,   use  Proposition \ref{Prop-22}  and integrate over $(0,T)$, we obtain
		\begin{align}\label{p01}
		\dis  -\int_0^T\langle \phi_t,\rho^n\rangle_{\V^\star,\V}\, dt + \int_0^T \mathcal{F}(\rho^n,\phi)dt
		=& \dis \int_{\omega_T}v^n\rho^n\, \phi \;\dq-\int_{\Om}\rho^n(T)\phi(T)\,\dx\notag\\
		&+\dis \int_\Omega \rho^0\,\phi(0) \dx.
		\end{align}	
Passing to the limit in \eqref{p01} while using \eqref{c21}, \eqref{c31} and \eqref{convB1}, we obtain	
		$$
		\dis \int_{\Om}\eta \phi(T) \,\dx -\int_0^T\langle \phi_t,\ro\rangle_{\V^\star,\V} \,\dx\, dt +\int_0^T \mathcal{F}(\ro,\phi)dt= \\
		\dis \int_{\omega_T}u\ro\, \phi \;\dq+\int_\Omega \rho^0\,\phi(0)\, \dx.
		$$
Using again Proposition \ref{Prop-22} we can deduce that
		\begin{equation}\label{SW1}
\begin{array}{rll}
	\dis \int_0^T\langle \ro_t+(-\Delta)^s\ro,\phi\rangle_{\V^\star,\V} \;dt
	&=&  \dis \int_{\omega_T}u\ro\, \phi \;\dq-\int_\Omega \rho^0\,\phi(0) \,\dx \\
&&+\dis \int_\Omega \rho^0\,\phi(0)\, \dx-\int_{\Om}\phi(T)(\eta-\ro(T))\dx,
		\end{array}
\end{equation}
for all $\phi \in W(0,T;\V)$. Since  $\ro$ is a solution of \eqref{modelc}, we get from \eqref{SW1} that
		$$
		\int_{\Om}\phi(T)(\eta-\ro(T))\, \dx=0, \quad
		\forall \phi \in W(0,T;\V).
		$$	
		Hence,  we can deduce that
		\begin{equation} \label{limp51}
		\eta=\ro(\cdot,T)\quad \text{ a.e. in}\quad \Om.
		\end{equation}
		Combining \eqref{c21}-\eqref{limp51} we obtain that, as $n\to\infty$,
		\begin{equation}\label{convroT}
		\ro^n(\cdot,T)\rightharpoonup \ro(\cdot,T)  \text{ weakly in }   L^2(\Omega).
		\end{equation}
Using \eqref{convroT}, \eqref{c11}, \eqref{ajout5} and the lower semi-continuity of $J$,  we can deduce that
			 $\dis J(u)\leq \liminf_{n\to \infty}J(v^n)=\inf_{v\in\U}J(v)$.
		This completes the proof.  		
	\end{proof}
We conclude this section with the following observation.	
	\begin{remark}
In Theorem \ref{existcontrol1} we only proved the existence of optimal solutions.  Uniqueness will necessitate additional assumptions.   As in the classical case under additional assumptions on the data $\rho^0,\rho^d\in L^\infty(\Omega)$ or on the parameter $\alpha$, then one can prove the uniqueness of local optimal solutions.
%
	\end{remark}
\section{Optimality conditions}\label{sec-5}
	In this section, we give the first and second order optimality conditions for the problem \eqref{optimal}, and hence, for  \eqref{modelc}-\eqref{pbopt}.
	\subsection{First order necessary optimality conditions}
The aim of this section is to derive the first order necessary optimality conditions and to characterize the optimal control. But before going further, we need some regularity results for the control-to-state operator.
Let us define the mapping
\begin{equation}\label{defG}
\left\{\begin{array}{l}
\mathcal G:W(0,T;\mathbb{ V})\times L^\infty(\omega_T)\to L^2((0,T);\V^\star)\times L^2(\Omega),\\
\mathcal{G}(\rho,v):= (\rho_t+(-\Delta)^s \ro-v\ro\chi_{\omega },\rho(0)-\rho^0).
\end{array}
\right.
	\end{equation}
Then, the state equation \eqref{modelc}  can be viewed as  $\mathcal{G}(\rho,v)=(0,0)$.
	\begin{lemma}\label{propd}
		The mapping $\mathcal{G}$ defined in \eqref{defG} is of class $\mathcal{C}^{\infty}$.
	\end{lemma}
\begin{proof}
	We write the first component $\mathcal{G}_1$ of $\mathcal{G}$  as
		$$\mathcal{G}_1(\rho,v)(\phi)=\mathcal{G}_{11}(\rho,v)(\phi)+\mathcal{G}_{12}(\rho,v)(\phi),\quad \forall \phi\in L^2((0,T);\V)$$
where
$$
\mathcal{G}_{11}(\rho,v)(\phi):=\int_0^T \left\langle\rho_t,\phi\right\rangle_{\V^\star,\V}\;dt+\int_0^T \mathcal{F}(\ro,\phi)dt$$
and
$$\mathcal{G}_{12}(\rho,v)(\phi):=\dis -\int_{\omega_T}v\ro\, \phi \;\dq.
$$
It is clear that  $\mathcal{G}_{11}$ is  linear and continuous  from $ W(0,T;\mathbb{ V})$ to $L^2((0,T);\V^\star)$ and $\mathcal{G}_{12}$ is  bilinear and continuous from  $ W(0,T;\mathbb{ V})\times L^\infty(\omega_T)$ to $L^2((0,T);\V^\star).$  Thus,  they are of class $\mathcal{C}^{\infty}$.
  The second component of $\mathcal{G}$ is clearly of class $\mathcal{C}^{\infty}$.
	\end{proof}
\begin{lemma}\label{lemmeG}
		The  mapping $G:L^\infty(\omega_T)\to W(0,T;\mathbb{ V}), u\mapsto \rho$ is of class $\mathcal{C}^{\infty}$.
	\end{lemma}
\begin{proof}
	Let $v\in L^\infty(\omega_T)$. It follows from Lemma \ref{propd} that $\mathcal{G}$ defined in \eqref{defG} is of class $\mathcal{C}^{\infty}$. Moreover,
		$$\partial_\rho\mathcal{G}(\rho,v)\varphi=\Big(\varphi_t+(-\Delta)^s \varphi-v\varphi\chi_{\omega },\varphi(0)\Big).$$
		 For  $\varphi^0\in L^2(\Omega)$ and $f\in L^2((0,T);\V^\star)$,   Corollary \ref{theoremexistence} shows that \eqref{model} with $\rho=\varphi$
		has a unique  solution $\varphi$ in $W(0,T;\mathbb{ V})$ which depends continuously on $\varphi^0$ and  $f$. Hence,
		$\dis \partial_\rho\mathcal{G}(\rho,v)$ defines an isomorphism from $W(0,T;\mathbb{ V})$ to $L^2((0,T);\V^\star)\times L^2(\Omega)$. Using the Implicit Function Theorem, we can deduce that $\mathcal{G}(\rho,v)=(0,0)$ has a unique solution $\rho=G(v)$ for any $v$ in $B^\infty_\varepsilon(u)$, where $ B^\infty_\varepsilon(u)$  denotes the open ball in  $L^\infty(Q)$
 of radius $\varepsilon$ centered at $u$,  solution of \eqref{pbopt}-\eqref{defuad}. Moreover,  the  operator $G:v\mapsto\rho$ is itself of class  $\mathcal{C}^{\infty}$. The proof is finished.
	\end{proof}

\begin{remark}
  Actually, we do not only have a unique solution $\rho=G(v)$ for any $v$ in a suitable neighborhood of $u$,  but we have a unique solution $\rho=G(v)$ of the state equation \eqref{modelc} for any given $v\in L^\infty(\omega_T)$ (see Corollary \ref{corexistence}). But this does not show that the optimal control found in Theorem \ref{existcontrol1} is unique.
\end{remark}

Next, we give some Lipschitz continuity results of the map $G$. The proof follows using similar arguments as in \cite{kenne2020,cb2021, cprg2021}.
	\begin{proposition}\label{prop2}
	Let  $v\in L^\infty(\omega_T)$ and $\rho^{0}\in L^\infty(\Omega)$. 
Then, there is a constant $C_1=C_1(N,s,\Omega)>0$  such that for all $v_1,v_2\in L^\infty(\omega_T)$,
\begin{equation}\label{estim20}
	\|G(v_1)-G(v_2)\|_{W(0,T;\V)}\leq C\|v_1-v_2\|_{L^2(\omega_T)}
\end{equation}
where $C=\left[\left( 2+C_1\|v_1\|_{\infty}\right)e^{\|v_1\|_{\infty}T}+1\right]
e^{\|v_2\|_{\infty}T}\|\rho^0\|_{L^\infty(\Omega)}.$	
	\end{proposition}
	\begin{proof}
		Let $v_1,v_2\in L^\infty(\omega_T).$
		Set $z:=\ro(v_1)-\ro(v_2)$, where $\ro(v_1),\ro(v_2)$ are  solutions of \eqref{modelc} with $v=v_1$ and $v=v_2$,
respectively. We have that  $z$ satisfies
		\begin{equation}\label{a1}
		\left\{
		\begin{array}{lllll}
		\dis  z_t+(-\Delta)^sz &=&(v_1z+(v_1-v_2)\ro(v_2))\chi_{\omega } &\hbox{in} & Q,\\
		\dis z&=&0 &\hbox{in}&\Sigma,\\
		z(\cdot,0)&=&0 &\hbox{in}& \Omega.
		\end{array}
		\right.
		\end{equation}
Since $v_1,v_2\in L^\infty(\omega_T)$ and $(v_1-v_2)\ro(v_2)\chi_{\omega }\in L^2(Q)$, using Corollary \ref{theoremexistence}, we have that there is a unique   $z\in W(0,T;\V)$ solution to  \eqref{a1}.   Thanks to \eqref{estimation22W0T}, we have that 		
\begin{equation}\label{rosny13}
\|z\|_{W(0,T;\V)}\leq
\left[\left( 2+C\|v_1\|_{\infty}\right)e^{\|v_1\|_{\infty}T}+1\right]\|(v_1-v_2)\ro(v_2)\|_{L^2(\omega_T)}.
 \end{equation}
Since $\rho^0\in L^\infty(\Omega)$,  we know that $\rho=\rho(v_2)$ solution of \eqref{modelc} satisfies
		\eqref{rhoinfini2},  that is,
		\begin{equation}\label{ajout4}
		\|\rho(v_2)\|_{L^\infty(\R^N\times (0,T))}\leq e^{\|v_2\|_{\infty}T}\|\rho^0\|_{L^\infty(\Omega)}.
		\end{equation}
Using the Cauchy-Schwarz inequality and \eqref{ajout4} in \eqref{rosny13}, we can deduce  \eqref{estim20}.
	\end{proof}
\begin{lemma}\label{diff}
		Let $G:L^\infty(\omega_T)\to W(0,T;\V)\cap L^\infty(\R^N\times(0,T)), v\mapsto \rho$ be the control-to-state operator,  where $\rho$ is the weak solution of \eqref{modelc}. Then, the directional derivative of $G$ in the direction $w\in L^\infty(\omega_T)$ is given by $G'(v)w=y,$
		where  $y\in W(0,T;\V)$  is the unique weak solution of
		\begin{equation}\label{diff1}
		\left\{
		\begin{array}{lllll}
		\dis  y_t+(-\Delta)^sy &=&(vy+w\ro)\chi_{\omega } &\hbox{in} & Q,\\
		\dis y&=&0 &\hbox{in}&\Sigma,\\
		y(\cdot,0)&=&0 &\hbox{in}& \Omega.
		\end{array}
		\right.
		\end{equation}	
Moreover, for every $v\in L^\infty(\omega_T)$, the linear mapping $w\mapsto G'(v)w$ can be extended to a linear continuous mapping from $L^2(\omega_T)\to W(0,T;\V)$.  In addition, there is a constant $C=C(N,s,\Omega)>0$ such that
\begin{equation}\label{estzvacter}
 \|y\|_{W(0,T;\V)} \leq
\left[\left( 2+C\|v\|_{\infty}\right)e^{\|v\|_{\infty}T}+1\right]
e^{\|v_2\|_{\infty}T}\|\rho^0\|_{L^\infty(\Omega)}\|w\|_{L^2(\omega_T)}.
\end{equation}	
	\end{lemma}
	\begin{proof}
	The proof is a consequence of the implicit function theorem.  In fact, to get \eqref{diff1} it is enough to differentiate the identity $\mathcal G(G(v),v)=(0,0)$.

 For the extension, it is sufficient to prove that for any $w\in L^2(\omega_T)$, the system \eqref{diff1} has a unique solution $y\in W(0,T;\V)$. This follows directly from Corollary \ref{theoremexistence}, with $f:=w\rho\in L^2(Q)$. Note that $f\in L^2(Q)$ because  $\rho\in L^\infty(\R^N\times(0,T))$.
	\end{proof}
	We have the following result.
	\begin{proposition}[\bf Fr\'echet differentiability of $J$]\label{diff4}
		Let $\ro$ be the solution of \eqref{modelc}.	Under the hypothesis of Lemma \ref{diff}, the functional $J:L^\infty(\omega_T)\to \R$ defined in \eqref{optimal} is  of class $C^\infty$  and for every $v,w\in L^\infty(\omega_T)$,
		\begin{equation}\label{diff5}
		J'(v)w=\int_{\Om}y(x,T)(\ro(x,T)-\ro^d(x))\,\dx+\alpha\int_{\omega_T}vw\ \dq,
		\end{equation}		
where $y\in W(0,T;\V)$ is the
		unique weak solution  of \eqref{diff1}.
	\end{proposition}
	\begin{proof}
		Observing that $\rho(v)=G(v)$ is a solution of \eqref{modelc}, 
		the first part of the proposition follows from 
		Lemma \ref{lemmeG}, since  $G$ has this property. The identity \eqref{diff5} follows from 	the chain rule.
		\end{proof}
Since the functional $J$ is non-convex, in general, we cannot expect a unique solution to the minimization problem \eqref{optimal}.  We introduce the following notion of local  solutions.
	\begin{definition}\label{defopt}
 We say that $u\in \U$ is an $L^\infty$-local solution
 of \eqref{optimal} if there exists $\varepsilon>0$ such that $J(u)\leq J(v)$ for every $v\in \U\cap B^\infty_\varepsilon(u)$.
	\end{definition}
The following result is crucial for the rest of the paper.
\begin{theorem}[\bf First order necessary optimality conditions]\label{theoSO}
	Let $\alpha>0,$  $v\in \mathcal{U}$ and $\rho^{0},\rho^d\in L^\infty(\Omega)$.  Let $u\in\mathcal U$ be an $L^\infty$-local minimum for  \eqref{optimal}. Then,
	\begin{equation}\label{ineq}
	J'(u)(v-u)\geq 0\;\;\;\text{for every}\;\;\; v\in \U.
	\end{equation}
	Moreover, there is a unique $\rho\in W(0,T;\V)$ and a unique $q\in W(0,T;\V)$  such that $\rho$ satisfies \eqref{modelc} with $v=u$,  $q$ satisfies	\eqref{ad21bis} with $v=u$
	and
	\begin{equation}\label{neccoptcond}
	\int_{\omega_T}(\alpha u+\ro q)(v-u)\ \dq\geq 0 \quad \forall v\in \mathcal{U}.
	\end{equation}
The condition \eqref{neccoptcond} is equivalent to the following:  for a.e. $(x,t)\in \omega_T$, we have that
\begin{equation}\label{contr1vac}
\dis \left\{\begin{array}{rlllll}
u(x,t)=m&\hbox{if}&\alpha u(x,t)+\ro(u(x,t))q(x,t)>0\\
u(x,t)\in [m,M] &\hbox{if}&\alpha u(x,t)+\ro(u(x,t))q(x,t)=0\\
u(x,t)=M&\hbox{if}&\alpha u(x,t)+\ro(u(x,t))q(x,t)<0.
\end{array}
\right.
\end{equation}
\end{theorem}
 \begin{remark} Note that \eqref{neccoptcond} or equivalently \eqref{contr1vac} can be rewritten as
	\begin{equation}\label{contr1}
	u=\min\left(\max\left(m,-\frac{q}{\alpha}\rho\right),M\right) \text{a.e.  in }\omega_T.
	\end{equation}	
\end{remark}
	\begin{proof}[\bf Proof of Theorem \ref{theoSO}]
	The proof of the first part of the theorem follows from the local case.
The equivalent between \eqref{neccoptcond} and \eqref{contr1vac} can be proved as in the local case contained in \cite[Chapter 4]{fredi2010}.
	\end{proof}
	\begin{remark}\label{contadj}
Using the change of variable $t\mapsto T-t$ and Lemma \ref{lemmeG}, one can show that the mapping $u\mapsto q $, solution of  \eqref{ad21bis} with $v=u$,  is also of class $\mathcal{C}^{\infty}$.
	\end{remark}
\begin{proposition}\label{propq}
Let  $u\in L^\infty(\omega_T)$ and $\rho^{0}, \rho^d\in L^\infty(\Omega)$.  
Let $q$ be the adjoint state solution to \eqref{ad21bis}  with $v=u$.
Then, there is a constant $C_1=C_1(N,s,\Omega)>0$ such that for all $u_1,u_2\in L^\infty(\omega_T)$,
\begin{equation}\label{estim20q}
		\dis \|q(u_1)-q(u_2)\|_{W(0,T;\V)}\leq  C\|u_1-u_2\|_{L^2(\omega_T)},
\end{equation}
where $C=e^{\|u_2\|_{\infty}T}\left[\left( 2+C_1\|u_1\|_{\infty}\right)e^{\|u_1\|_{\infty}T}+1\right]
\left(\|\rho^0\|_{L^\infty(\Omega)}+\|\ro^d\|_{L^\infty(\Omega)}\right)$.
\end{proposition}
\begin{proof}
We proceed as in  the proof of Proposition \ref{prop2}. Let $u_1,u_2\in L^\infty(\omega_T).$
		Set $z:=q(u_1)-q(u_2)$, where $q(u_1),q(u_2)\in W(0,T;\V)$ are solutions of \eqref{ad21bis} with $v=u_1$ and $v=u_2$,
respectively. We have that  $z$ satisfies
		\begin{equation}\label{a1-1}
		\left\{
		\begin{array}{lllll}
		\dis  z_t+(-\Delta)^sz &=&(u_1z+(u_1-u_2)q(u_2))\chi_{\omega } &\hbox{in} & Q,\\
		\dis z&=&0 &\hbox{in}&\Sigma,\\
		z(\cdot,0)&=&0 &\hbox{in}& \Omega.
		\end{array}
		\right.
		\end{equation}
Using Corollary \ref{theoremexistence}, we have that there exists a unique   $z\in W(0,T;\V)$ solution to  \eqref{a1-1}.   Thanks to  \eqref{estimation22W0T}, we have that there is a constant $C>0$ such that
\begin{equation}\label{rosny13bis}
		\dis \|z\|_{W(0,T;\V)}\leq
  \left[\left( 2+C\|u_1\|_{\infty}\right)e^{\|u_1\|_{\infty}T}+1\right]\|(u_1-u_2)q(u_2)\|_{L^2(\omega_T)}.
		\end{equation}	
Since $\rho(\cdot,T)-\rho^d\in L^\infty(\Omega)$ where $\rho$ is the weak solution of \eqref{modelc},  we know that $q=q(u_2)$ solution of \eqref{ad21bis} satisfies \eqref{qW0T},  that is,
		\begin{equation}\label{ajout4bis}
		\|q(u_2)\|_{L^\infty(\R^N\times (0,T))}\leq e^{\|u_2\|_{\infty}T}\left(\|\rho^0\|_{L^\infty(\Omega)}+\|\ro^d\|_{L^\infty(\Omega)}\right).
		\end{equation}
Using the Cauchy-Schwarz inequality and \eqref{ajout4bis} in \eqref{rosny13bis}, we can deduce \eqref{estim20q}.
\end{proof}

We conclude this section with the following result.

\begin{lemma}\label{cont1}
	For every $u\in L^\infty(\omega_T)$,  the linear mapping $v\mapsto J'(u)v$ can be extended to a linear continuous mapping $J'(u):L^2(\omega_T)\to \R$ given by \eqref{diff5}.
\end{lemma}	
\begin{proof} Let $u\in L^\infty(\omega_T)$ and $v\in L^2(\omega_T)$. From \eqref{diff5}, we have that
	$$
	J'(u)v=\int_{\omega_T}(\alpha u+\rho q)v\, \dq,
	$$
	where $\rho$ and $q$ are  solutions of \eqref{modelc} and \eqref{ad21bis} with $v=u$, respectively.
	Using \eqref{rhoinfini2} and \eqref{qW0T},  we have that there is a constant $C>0$ independent of $v$ such that
	$$
	\begin{array}{lll}
	\dis |J'(u)v|\leq C\|v\|_{L^2(\omega_T)}.
	\end{array}
	$$
Thus, the mapping $v\mapsto J'(u)v$ is  linear and continuous  on $L^2(\omega_T)$.
\end{proof}

	\subsection{Second order necessary and sufficient optimality conditions}\label{sufficientoptcond}
Note that
the cost functional $J$ associated to the optimization problem \eqref{optimal} is non-convex and the first order optimality conditions given in Theorem \ref{theoSO} are necessary but not sufficient for optimality. The sufficiency requires the use of second order optimality conditions, which is the aim of this section. To proceed, we need the following result.
	\begin{lemma}[\bf The mapping $G$ is of class $C^\infty$]\label{diff2} Let $\alpha>0,$  $v\in \mathcal{U}$ and $\rho^{0},\rho^d\in L^\infty(\Omega)$.  Let $u$ be an $L^\infty$-local minimum for the  problem \eqref{optimal}. Then,  the control-to-state mapping $G:u\mapsto \rho$  is of class $C^\infty$
	from $L^\infty(\omega_T)$ into $W(0,T;\V)$.   Moreover,
$G''(u)[w,h]=z,$
		where $w,h\in L^\infty(\omega_T)$ and $z\in W(0,T;\V)$ is the unique weak solution of
		\begin{equation}\label{diff3}
		\left\{
		\begin{array}{lllll}
		\dis  z_t+(-\Delta)^sz &=&(uz+hG'(u)w+wG'(u)h)\chi_{\omega } &\hbox{in} & Q,\\
		\dis z&=&0 &\hbox{in}&\Sigma,\\
		z(\cdot,0)&=&0 &\hbox{in}& \Omega.
		\end{array}
		\right.
		\end{equation}		
	\end{lemma}
\begin{proof}
The first part of the lemma is a direct consequence of Lemma \ref{lemmeG} and \eqref{diff3} is a consequence of the Implicit Function Theorem.  We omit the details for brevity.
\end{proof}
	\begin{lemma}[\bf The mapping $J$ is of class $C^\infty$]\label{diff44}
	Let $\alpha>0$   and $\rho^{0},\rho^d\in L^\infty(\Omega)$.  Let $u\in\mathcal U$ be an $L^\infty$-local minimum for the minimization problem \eqref{optimal} associated to the state $\ro:=G(u)$ solution of \eqref{modelc}. Let $q$ be the solution of \eqref{ad21bis} with $v=u$.	Then, the functional $J: L^\infty(\omega_T) \to \R$ is of class $C^\infty$ and for every $w, h\in L^\infty(\omega_T)$, we have that
		\begin{align}\label{diff6}
		J''(u)[w,h]=&\int_{\omega_T}[hG'(u)w+wG'(u)h]q\, \dq+\int_{\Om}(G'(u)w)(x,T)(G'(u)h)(x,T)\, \dx\notag\\
		&+\alpha\int_{\omega_T}hw\, \dq.
		\end{align}	
		Moreover, 	the bilinear mapping $(w,h)\mapsto J''(u)[w,h]$ can be extended to a bilinear continuous mapping $J''(u):L^2(\omega_T)\times L^2(\omega_T)\to \R$ given by \eqref{diff6}.
	\end{lemma}
	
	\begin{proof}
The identity \eqref{diff6} and the fact that $J$ is of class $C^\infty$  follow from Lemmas \ref{lemmeG},  \ref{diff2},  a straightforward computation and the Lebesgue dominated convergence theorem.  For completeness we give the details of the proof of 
the last part.   Let $w,h\in L^2(\omega_T)$.
Recall that $G'(u)w$ is a solution of \eqref{diff1}.
	Using Cauchy-Schwarz's inequality, \eqref{estzvacter}  and \eqref{qW0T}, we obtain that there is a constant $C=C(N,s, \Omega)>0$ such that
\begin{align*}
	\dis |J''(u)[w,h]|	\leq& \|q\|_{L^\infty(\R^N\times (0,T))}\|h\|_{L^2(\omega_T)}\|G'(u)w\|_{L^2(Q)}\\
&+\|q\|_{L^\infty(\R^N\times (0,T))}\|w\|_{L^2(\omega_T)}\|G'(u)h\|_{L^2(Q)}\\
	&+ \|(G'(u)w)(T)\|_{L^2(\Omega)}\|(G'(u)h)(T)\|_{L^2(\Omega)}\\
&+\alpha \|w\|_{L^2(\omega_T)}\|h\|_{L^2(\omega_T)}\\
	\leq& \|h\|_{L^2(\omega_T)}\|w\|_{L^2(\omega_T)}\left(2 C_2\|q\|_{L^\infty(\R^N\times (0,T))}+C_2^2+\alpha\right)
\\
		\leq&C_3\|w\|_{L^2(\omega_T)}\|h\|_{L^2(\omega_T)}
	\end{align*}
where $$C_3=\left( C_2e^{\|u\|_{\infty}T}\left(\|\rho^0\|_{L^\infty(\Omega)}+\|\ro^d\|_{L^\infty(\Omega)}\right)
+C_2^2+\alpha\right)$$ with
$C_2=\left[\left( 2+C\|u\|_{\infty}\right)e^{\|u\|_{\infty}T}+1\right]
e^{\|u\|_{\infty}T}\|\rho^0\|_{L^\infty(\Omega)}$.
	Hence,  the mapping $(w,h)\mapsto J''(u)[w,h]$ is bilinear and continuous on $L^2(\omega_T)\times L^2(\omega_T)$.
	\end{proof}
	
Next, we introduce some  concepts retrieved from \cite{fredi2010}.
For a given $\tau\geq 0$, we define the set of strongly active constraints  $A_\tau(u)$  by
	\begin{equation*}
		A_\tau(u):=\{(t,x)\in \omega_T : |\alpha u(x,t)+\rho(u)q(x,t)|>\tau\}.
		\end{equation*}

The $\tau$-critical  set  associated to a control $u$ (see e.g. \cite{fredi2010}) is defined by
		\begin{equation}\label{ccone}
		\mathcal C_\tau (u)= \{v\in L^\infty(\omega_T) : v \;\text{fulfills}\; \eqref{eq1}\},
		\end{equation}
that is, for a.e $(t,x)\in Q$,  we have that
		\begin{equation}\label{eq1}
		\left\{
		\begin{array}{lllll}
		\dis  v(x,t)\geq 0 &if&  u(t,x)=m \hbox{ and } (t,x)\notin A_\tau(u),\\
		\dis v(x,t)\leq 0&if& u(t,x)=M \hbox{ and } (t,x)\notin A_\tau(u),\\
		v(x,t)=0&if&(t,x)\in A_\tau(u).
		\end{array}
		\right.
		\end{equation}
	
	We have the following observation as in \cite{casas2015a}.
	\begin{remark}
		We notice that from the differentiability in $L^\infty(\omega_T)$, the cone defined in \eqref{ccone} is a subset of $L^\infty(\omega_T)$.  Since $L^\infty(\omega_T)$ is dense in $L^2(\omega_T)$ and from Lemma \ref{diff44}, the quadratic form $ J''(u)$ is continuous on $L^2(\omega_T)$, we have that the second order conditions based on the critical cone in $L^\infty(\omega_T)$ can be transfered to the extended cone in $L^2(\omega_T)$. This cone is convex and closed in $L^2(\omega_T)$.
	\end{remark}
	
In the rest of the paper, we will adopt the notation $J''(u)v^2:=J''(u)[v,v]$.
\begin{theorem}[\bf Second order necessary optimality conditions]\label{theo-514}
	Let $u\in \mathcal{U}$ be an $L^\infty$-local solution of \eqref{optimal}. Then,  $J''(u)v^2\geq 0$ for all $v\in \mathcal C_0 (u)$.
\end{theorem}

\begin{proof}
The proof follows as in  \cite[pp. 246]{fredi2010}.  We omit the details for brevity.
\end{proof}

\begin{lemma}\label{hypothesis}
	Let $u\in \mathcal{U}$ be an $L^\infty$ local-control satisfying  \eqref{ineq}. Then, the following assertions hold:
	\begin{enumerate}
		\item The functional $J:L^\infty(\omega_T)\to \R$ is of class $\mathcal{C}^\infty$. Furthermore, there exist continuous extensions
		\begin{equation}\label{h1}
		J'(u)\in \mathcal{L}(L^2(\omega_T),\R)\;\;\;\;\text{and}\;\;\;J''(u)\in \mathcal{B}(L^2(\omega_T),\R).
		\end{equation}
		\item For any sequence $\left\{(u_k,v_k)\right\}_{k=1}^{\infty}\subset \mathcal{U}\times L^2(\omega_T)$ with $\|u_k-u\|_{L^2(\omega_T)}\to 0$ and $v_k \rightharpoonup v$ weakly in $L^2(\omega_T)$,  as $k\to\infty$, we have that
		\begin{equation}\label{h2}
		\dis	J'(u)v=\lim_{k\to \infty} J'(u_k)v_k,
		\end{equation}
and		
		\begin{equation}\label{h3}
		\dis	J''(u)v^2\leq \liminf_{k\to \infty} J''(u_k)v_k^2.
		\end{equation}
If  $v=0$,  then
		\begin{equation}\label{h4}
\alpha \liminf_{k\to \infty}\|v_k\|^2_{L^2(\omega_T)}\leq\liminf_{k\to \infty} J''(u_k)v_k^2.
		\end{equation}
	\end{enumerate}
\end{lemma}

\begin{proof}
(a) This part follows from Proposition \ref{diff4},  Lemmas \ref{cont1} and \ref{diff44}.

(b)  We proceed in three steps.
	
 \textbf{Step 1.} We show \eqref{h2}.  Using \eqref{estim20} in Proposition \ref{prop2},  we get that $\dis G(u_k)\to G(u)$ in $W(0,T;\V)$,  as $k\to \infty$. From Proposition \ref{propq}, we also obtain that the adjoint state $q$ solution of \eqref{ad21bis} satisfies $\dis q(u_k)\to q(u)$ in $W(0,T;\V)$,  as $k\to \infty$.

 We claim that
	\begin{equation}\label{conv}
	G(u_k)q(u_k)\to G(u)q(u) \;\;\text{in}\;\;L^2(Q), \;\;\text{as }\;\; k\to \infty.
	\end{equation}
Indeed, notice that $G(u_k)$ and $G(u)$ are solutions to \eqref{modelc} with $v=u_k$ and $v=u$, respectively. Therefore, from \eqref{rhoinfini2}, we can deduce that $G(u_k)$ is bounded in $L^\infty(Q)$ and $G(u)\in L^\infty(Q)$. Thus,  $G(u_k)q(u_k),G(u)q(u)\in L^2(Q)$ and we have that
\begin{align}\label{mm}
\|G(u_k)q(u_k)-G(u)q(u)\|_{L^2(Q)}\leq& \|q(u_k)-q(u)\|_{L^2(Q)}\|G(u_k)\|_{L^\infty(Q)}\notag\\
&+\|G(u_k)-G(u)\|_{L^2(Q)}\|q(u)\|_{L^\infty(Q)} .
\end{align}	
Taking the limit,  as $k\to\infty$,  of \eqref{mm} we obtain the claim \eqref{conv}.
From \eqref{conv} we can deduce that $\alpha u_k+G(u_k)q(u_k)\to \alpha u+G(u)q(u)$ in $L^2(Q)$,  as $k\to \infty$.  Using the expression of $J'$ given in \eqref{diff5}, we have that
	\begin{align*}
\lim_{k\to \infty} J'(u_k)v_k=&\lim_{k\to \infty}\int_{\omega_T}(\alpha u_k+G(u_k)q(u_k))v_k\, \dq \\
=&\int_{\omega_T}(\alpha u+G(u)q(u))v\, \dq= {J}'(u)v.
	\end{align*}
We have shown \eqref{h2}.
	
 \textbf{Step 2.} We show \eqref{h3}. We write
	\begin{align}\label{d}
	J''(u_k)v_k^2=&2\int_{\omega_T}v_k(G'(u_k)v_k)q(u_k)\dq+\int_{\Om}|(G'(u_k)v_k)(T)|^2\, \dx\notag\\
	&+\alpha\int_{\omega_T}|v_k|^2\, \dq,
	\end{align}
where	$G'(u_k)v_k$ is the unique weak solution of \eqref{diff1} with $u=u_k$ and $w=v_k$. Since $G(u_k)$ and $v_k$ are bounded in $L^{\infty}(Q)$ and in $L^2(\omega_T)$,  respectively,  it follows from  \eqref{estzvacter} that $G'(u_k)v_k$ and $(G'(u_k)v_k)(\cdot,T)$ are  bounded in $W(0,T;\V)$ and in $L^2(\Omega)$, respectively. Thus, up to a subsequence if necessary,  as $k\to \infty$,  $G'(u_k)v_k$ converges weakly to $G'(u)v$ in $L^2((0,T);\V)$ and $(G'(u_k)v_k)(\cdot,T)$ converges weakly to $(G'(u)v)(\cdot,T)$ in $L^2(\Omega)$. Now, thanks to Theorem \ref{thmcompact},  the embedding $W(0,T;\V)\hookrightarrow L^2(Q)$ is compact. Hence, we obtain that $G'(u_k)v_k$ converges strongly to $G'(u)v$  in $L^2(Q)$,  as $k\to \infty$.
Also,  since $\|u_k-u\|_{L^2(\omega_T)}\to 0$  as $k\to \infty$, we have from Proposition \ref{propq} that  $q(u_k)\to q(u)$ in $L^2((0,T);\V)$,  as $k\to \infty$,  and $q(u_k)$ is bounded in $L^\infty(Q)$ (see Lemma \ref{corominmax2}). Therefore, $(G'(u_k)v_k) q(u_k)$ converges strongly to $(G'(u)v)q(u)$  in $L^2(Q)$,  as $k\to \infty$.  Taking the limit,  as $k\to \infty$,  in \eqref{d} and using the lower-semi continuity of the $L^2$-norm, we can deduce that
	\begin{align*}
		\dis \lim_{k\to \infty}  J''(u_k)v_k^2\ge& 2\lim_{k\to \infty}\int_{\omega_T}v_k(G'(u_k)v_k)q(u_k)\dq\\
		&+\liminf_{k\to \infty}\left[\int_{\Om}|(G'(u_k)v_k)(T)|^2\, \dx+\alpha \int_{\omega_T}|v_k|^2\, \dq\right]\\
		\geq& 2\int_{\omega_T}v(G'(u)v)q(u)\dq+\int_{\Om}|(G'(u)v)(T)|^2\, \dx+\alpha \int_{\omega_T}|v|^2\, \dq\\
		=& J''(u)v^2.
	\end{align*}
We have shown \eqref{h3}.
	
 \textbf{Step 3.}
	If $v=0$, then in \eqref{d} the first and second terms tend to $0$.  Hence,
	$$\alpha  \liminf_{k\to \infty}\|v_k\|^2_{L^2(\omega_T)}\leq   \lim_{k\to \infty}  J''(u_k)v_k^2.$$ The proof is finished.
\end{proof}

Now, from Lemma \ref{hypothesis}, we obtain that the assumptions in \cite[Theorem 2.3 ]{casas2012} are fulfilled.  Hence, we obtain the following result which provides second-order sufficient conditions for locally optimal solutions.

\begin{theorem}[\bf Second order sufficient optimality conditions]\label{Quadratic growth}
Let $\rho^0,\rho^d\in L^\infty(\Omega)$.
	Let $u\in \mathcal{U}$ be a control satisfying \eqref{neccoptcond} and
\begin{align}\label{CD}
J''(u)v^2>0 \quad \forall v\in \mathcal C_0 (u)\setminus \{0\}.
\end{align}
Then, there are two constants $\gamma>0$ and $\beta>0$ such that
	\begin{equation}\label{cdt2}
	J(v)\geq J(u)+\frac{\beta}{2}\|v-u\|^2_{L^2(\omega_T)}\quad \forall v\in \mathcal{U}\cap B^2_{\gamma}(u),
	\end{equation}
	where $B^2_{\gamma}(u)$ is the open ball in $L^2(\omega_T)$ with center $u$ and radius $\gamma$.
\end{theorem}


\begin{remark}
In this work, we proved that the cost functional $J$ is of class $\mathcal{C}^{\infty}$ only in $L^\infty(\omega_T)$. Hence,  the so-called two-norms discrepancy (see e.g. \cite{casas2012}) occurs.  The two-norms discrepancy occurs when the functional $J$ is twice differentiable with respect to one norm (the $L^\infty$-norm in our case), but the lower estimates for the second derivative $J''$ holds in a weaker norm (the $L^2$-norm in our case) in which $J$ is not twice differentiable. For more details on the topic, we refer to \cite{antil2020, ioffe1979, otarola2022} and the references therein.
\end{remark}

\subsection*{Data availability}	 Data sharing not applicable to this article as no datasets were generated or analyzed during the current study.
	\section*{Declarations}
\subsection*{Conflict of interests}
 The authors have no competing interests to declare that are relevant to the content of this article.
	
	\bibliographystyle{plain}
	\bibliography{references}

\end{document}